\documentclass[11pt]{amsart}

\usepackage{epsfig}

\usepackage{graphicx}
\usepackage{amscd}
\usepackage{amsmath}
\usepackage{amsfonts}
\usepackage{yfonts}
\usepackage{psfrag}
\usepackage{amssymb}
\usepackage{verbatim}
\usepackage{comment}

\newtheorem{theorem}{Theorem}[section]
\theoremstyle{plain}

\newtheorem{lemma}[theorem]{Lemma}

\def\diam{{\rm diam}}

\def\th{\theta}

\def\wtil{\widetilde}

\newcommand{\lam}{\lambda}

\newcommand{\gam}{\gamma}

\newcommand{\R}{{\mathbb R}}
\newcommand{\Q}{{\mathbb Q}}
\newcommand{\Z}{{\mathbb Z}}
\newcommand{\C}{{\mathbb C}}
\def\Compl{{\mathbb C}}

\newcommand{\Nat}{{\mathbb N}}

\newcommand{\m}{{\bf m}}

\def\Ck{{\mathcal C}}
\def\Tk{{\mathcal T}}
\def\Vk{{\mathcal V}}

\newcommand{\eps}{{\varepsilon}}
\newcommand{\es}{\emptyset}

\def\ov{\overline}
\newcommand{\const}{{\rm const}}
\def\Span{{\rm Span}}

\begin{document}

\title[Expansion maps for self-affine tilings]{On the 
characterization of expansion maps for self-affine tilings}

\author{Richard Kenyon}
\address{Richard Kenyon\\ Department of Mathematics\\ Brown University\\ Providence, RI 02912}
\author{Boris Solomyak }
\address{Boris Solomyak, Box 354350, Department of Mathematics,
University of Washington, Seattle WA 98195}
\email{solomyak@math.washington.edu}

\begin{abstract}
We consider self-affine tilings in $\R^n$ with expansion matrix $\phi$ 
and address the question which matrices $\phi$
can arise this way. In one dimension, $\lam$ is an expansion factor of
a self-affine tiling if and only if $|\lam|$ is a Perron number,
by a result of Lind. In two dimensions, when $\phi$ is
a similarity, we can speak of a complex expansion factor, and there is
an analogous necessary condition,
due to Thurston: if a complex
$\lam$ is an expansion factor of a self-similar tiling, then 
it is a complex Perron number. We establish a necessary condition for
$\phi$ to be an expansion matrix for any $n$, assuming only that $\phi$ is
diagonalizable over $\C$. We conjecture that this condition on
$\phi$ is also sufficient for the existence of a self-affine tiling.
\end{abstract}

\date{\today}

\thanks{
The research of R. K.  was supported in part by NSERC.
The research of B. S.  was supported in part by NSF 
}

\maketitle

\thispagestyle{empty}

\section{Introduction}

Self-affine tilings arise in many different contexts, notably in dynamics
(Markov partitions for hyperbolic maps \cite{Sinai,KV,Prag}), 
logic (aperiodic tilings \cite{Penrose}), number theory
(radix representations \cite{Rau,LW}), physics  (quasicrystals
\cite{BomTay}), 
ergodic theory \cite{Solomyak}, and hyperbolic groups
\cite{Wordprocessing}.  See \cite{BBL,Robi} for recent surveys with a
large bibliography.

A {\bf self-affine tiling} (SAT) $\Tk= \{T_i\}_{i\in I}$ 
of $\R^n$ is a covering of $\R^n$ with sets (tiles) $T_i$ 
satisfying the following properties:
\begin{enumerate}
\item Each tile $T_i$ is the closure of its interior.
\item Interiors of tiles do not overlap.
\item There are a finite number of tile types up to translation.
\item The tiling is {\bf repetitive} and has {\bf finitely many local
configurations} (see the next section for definitions).
\item There is an expanding linear map $\phi:\R^n\to\R^n$ mapping tiles over tiles:
the image of a tile $T_i$ is a union of tiles of $\Tk$, and two tiles of the same
type have images which are translation-equivalent patches of tiles. 
\end{enumerate}

\begin{figure}[htbp]
\center{\scalebox{1.0}{\includegraphics{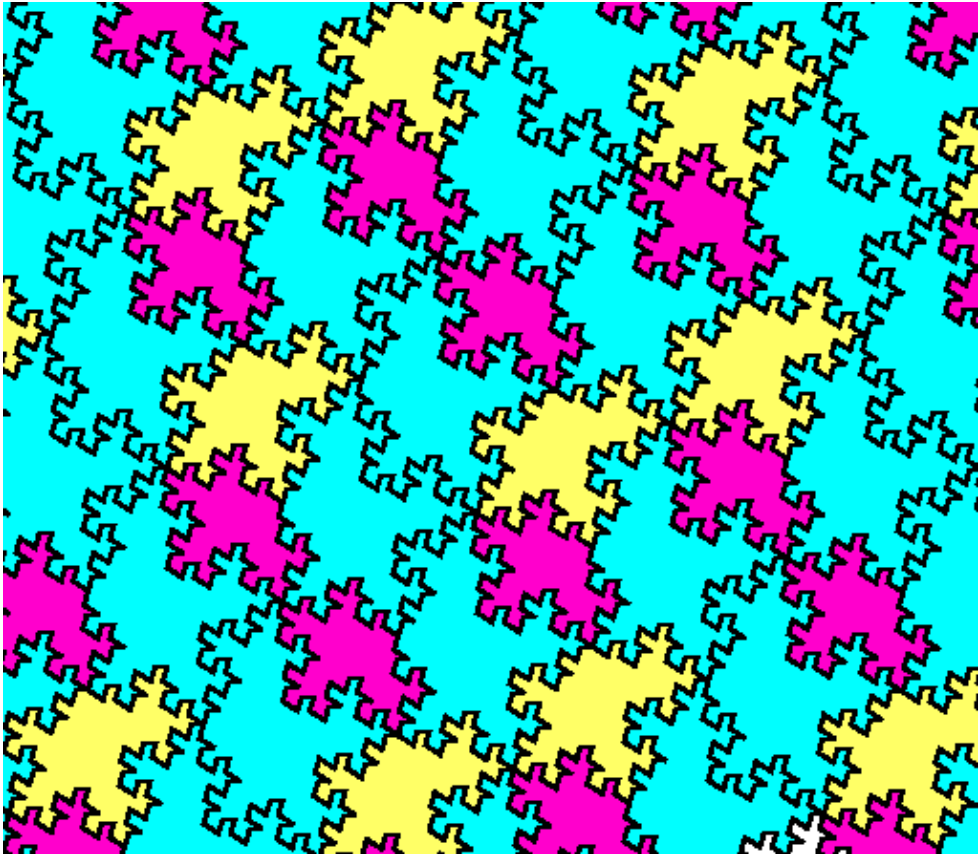}}}
\caption{\label{SAT1}A self-affine tiling in the plane with expansion $\phi(z)=\lambda z$ where $\lambda$ is the complex root of $x^3+x+1=0$. Here there are three tile types, all similar. The smallest
scales to the medium which scales to the large; the large
subdivides into a small and a large. One can construct this
tiling using the method of \cite[Sec.6]{Kenyon.construction},
as follows. To a reduced word in the free
group on three letters $F(a,b,c)$ associate a polygonal
path in $\C$ by sending $a^{\pm 1}$ to $\pm 1$,
$b^{\pm1}$ to $ 
\pm \lambda$, $c^{\pm 1}$ to $\pm \lambda^2.$ 
Let $\psi$ be the endomorphism of 
$F(a,b,c)$ defined by $\psi(a)=b,\psi(b)=c,\psi(c)=a^{-1}b^{-1}.$ Consider the three commutators $[a,b]=aba^{-1}b^{-1}, [b,c],$ and $[a,c]$; they represent three
closed paths.
Then 
$\lim_{n\to\infty}\lambda^{-n}\psi^n([a,c])$ is the boundary of
the smallest tile; the other tiles boundaries are
$\lim_{n\to\infty}\lambda^{-n}\psi^n([a,b])$ and 
$\lim_{n\to\infty}\lambda^{-n}\psi^n([b,c]).$
The subdivision rule comes from the identities 
$\psi([a,c])=a^{-1}[a,b]a,$
$\psi[a,b]=[b,c]$ and $\psi[b,c]=[c,a^{-1}b^{-1}]=
(a^{-1}[a,c]a)(a^{-1}b^{-1}[b,c]ba).$
}
\end{figure}

\begin{figure}
\center{\scalebox{.8}{\includegraphics{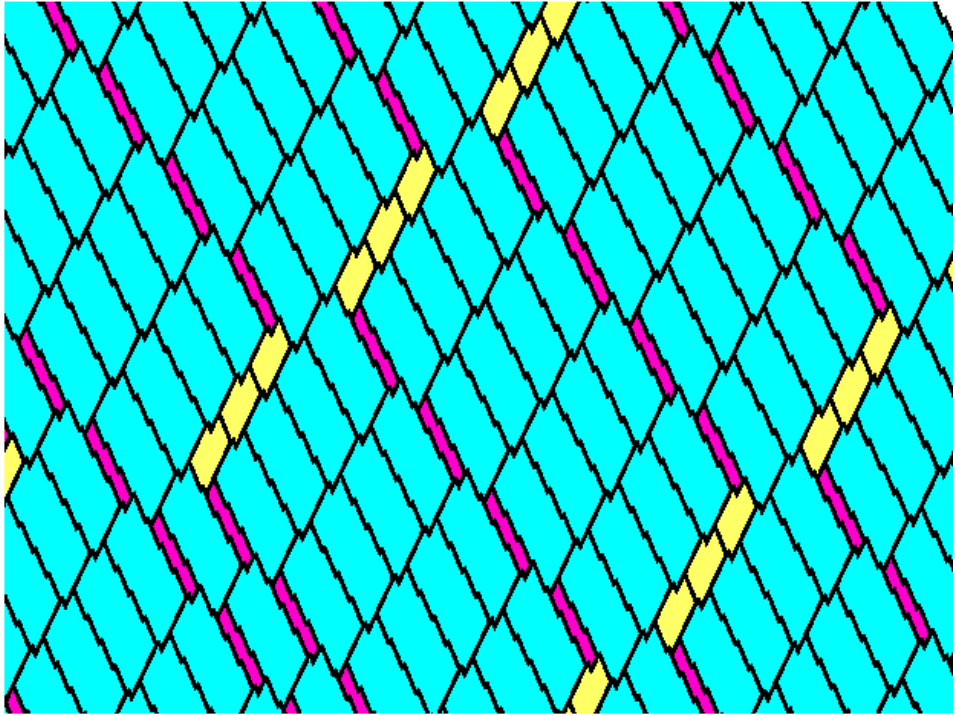}}}
\caption{\label{SAT2}A self-affine tiling in the plane with diagonal expansion
matrix ${\rm Diag}[x_1, x_2]$ where
$x_1\approx 2.19869, \ x_2 \approx - 1.91223$ are roots of
$x^3-x^2-4x+3=0$.}
\end{figure}

The simplest example is the periodic tiling with unit cubes
and expansion mapping $\phi(x)=2x$. However 
typically SATs are nonperiodic and have tiles with fractal 
boundaries. See Figures \ref{SAT1} and \ref{SAT2} for examples in $\R^2$.

Lind \cite{Lind} (using different language) 
gives a characterization of expansion factors of self-affine tilings in one dimension: $\lambda$ is the expansion of an SAT of $\R$ if and only if
$|\lambda|$ 
is a {\bf Perron number}, that is, a real algebraic integer which is strictly
larger in modulus than all of its Galois conjugates. 

A self-affine tiling is {\bf self-similar} if $\phi$ is a similarity (a homothety followed
by a rotation).
Thurston \cite{Thurston} showed that the expansion factor $\lambda\in\C$ of
a self-similar tiling of $\R^2$ is a {\bf complex Perron number}, that is,
an algebraic integer which is strictly larger in modulus than its Galois conjugates
except for its complex conjugate. 
In \cite{Kenyon.construction}, a construction
of a self-similar tiling for every complex Perron number is given;
unfortunately, the proof as written in subsection 4.5
of \cite{Kenyon.construction} is incomplete. 
A version of the construction does yield a tiling with  expansion
$\lambda^k$ for $k$ sufficiently large, and we hope that  it can be
modified to get a tiling with expansion $\lam$, completing the characterization.
This gap does not affect the construction in 
section 6 of \cite{Kenyon.construction}
which uses free group endomorphisms;
however, the latter does not cover all the 
complex Perron numbers. See also \cite{FIR} for a related construction.

In the current paper we study SATs of $\R^n$ with expansion matrix $\phi$ which
is diagonalizable over $\C$. We show that if $\phi$ is the expansion matrix for an SAT
then eigenvalues of $\phi$ are algebraic integers, and for every eigenvalue
$\gamma$, all Galois conjugates of $\gamma$ which have modulus
$\geq|\gamma|$ have multiplicity (among eigenvalues of $\phi$)
at least as large as that of $\gamma$, see Theorem \ref{th-main} below.

An alternative description of this criterion is that there is an integer matrix
$M$ acting on $\R^N$ for some $N\geq n$, which has an invariant real subspace $W$
of dimension $n$, on which it has strictly larger growth (that is,
strictly larger determinant, in absolute value) than for any other
$n$-dimensional invariant subspace, and $M$ restricted to $W$ is linearly conjugate to $\phi$.

The converse to our result is open: does there exist, for every
linear map $\phi$ satisfying the above conditions, an SAT
with expansion $\phi$? We conjecture that the answer is yes.

In  Figure 2 we show an example of a self-affine (non-self-similar)
SAT in the plane.
The subdivision rule
is indicated in Figure 3.

\begin{figure}
\scalebox{.5}{\includegraphics{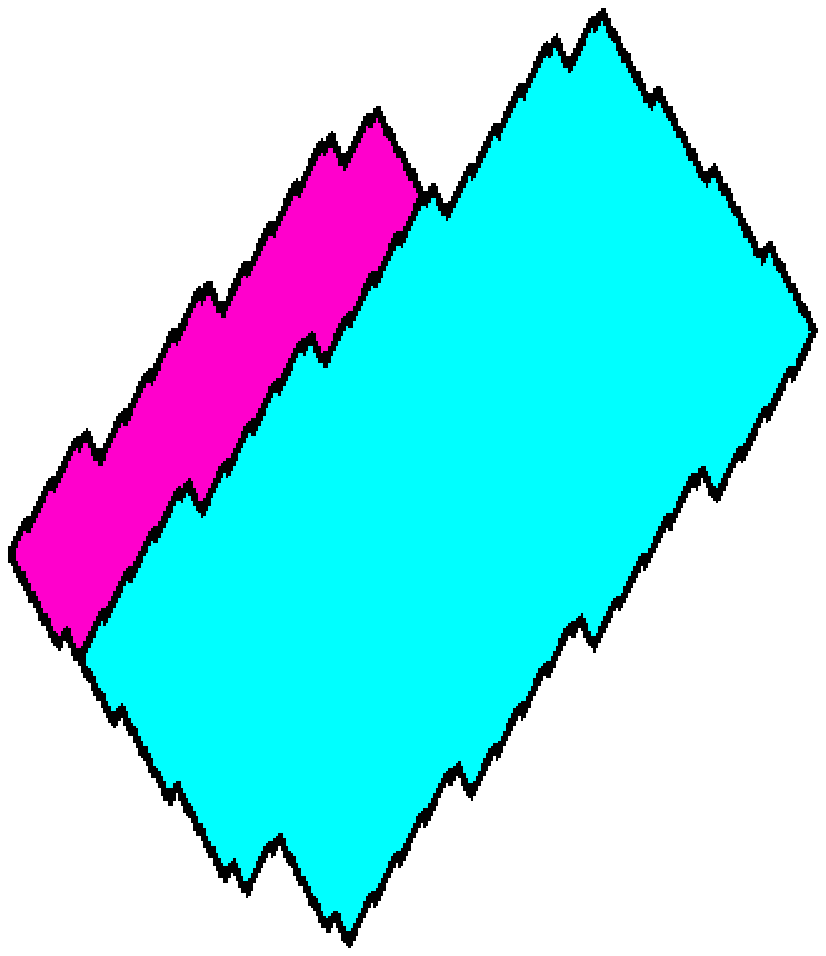}}
\scalebox{.5}{\includegraphics{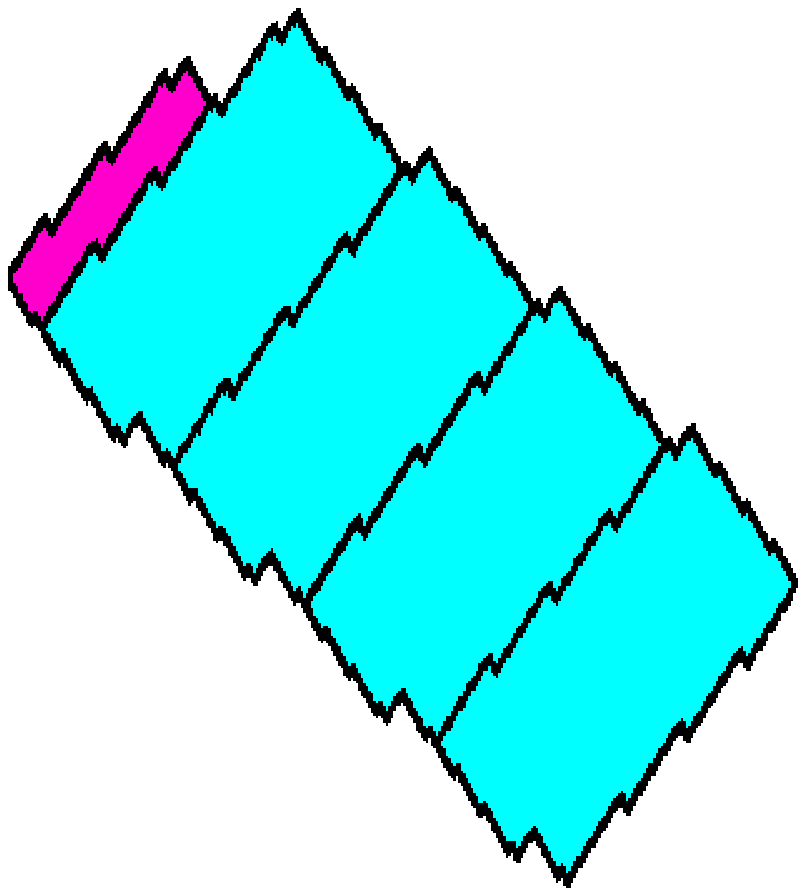}}
\scalebox{.5}{\includegraphics{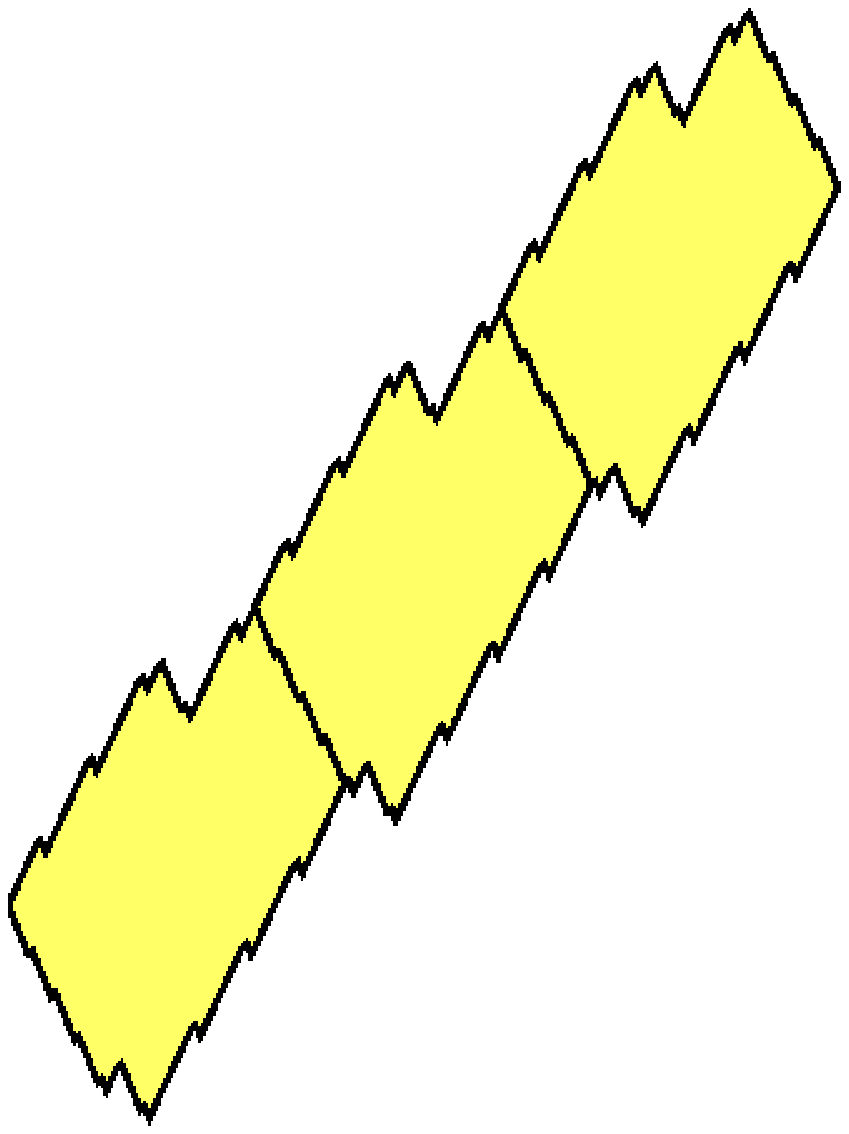}}
\caption{Subdivision rule: $1 \to \{3,2\}$, $2 \to \{3,2,2,2,2\}$,
$3 \to \{1,1,1\}$. The construction is similar to the previous example
but with $a,b,c$ corresponding to vectors
$(1,1),(x_1-1,x_2-1),(x_1^2-x_1,x_2^2-x_2)$ in $\R^2$, endomorphsm $\psi(a)=ab,\psi(b)=c,\psi(c)=ab^4$ and
tiles $[b,a],[b,c],[a,c].$}
\end{figure}

Our methods do not at present extend to the non-diagonalizable case.
However, we conjecture that the second description above holds in general, that is,
without the constraint of diagonalizability, $\phi$ is the expansion of an SAT if and only
if there is an integer matrix $M$ acting on $\R^N$ for some $N\geq n$,
which has an invariant real subspace $W$
of dimension $n$, on which it has strictly larger growth (determinant) than for any other
$n$-dimensional invariant subspace, and $M$ restricted to $W$ is linearly conjugate to $\phi$.
For example, we conjecture that there is no SAT in $\R^3$ with expansion
$$\left(\begin{matrix}3+\sqrt{2}&1&0\\0&3+\sqrt{2}&0\\0&0&3-\sqrt{2}\end{matrix}\right)$$
although it is easy to construct one with expansion
$$\left(\begin{matrix}3+\sqrt{2}&0&0\\0&3+\sqrt{2}&0\\0&0&3-\sqrt{2}\end{matrix}\right)$$

\section{Preliminaries}

We say that a tiling $\Tk=\{T_i\}_{i\in I}$ has a finite number of tile
types up to translation, if there is an equivalence relation $\sim$ on the tiles $T_i$
with a finite number of equivalence classes and $T_i\sim T_j$ implies
that $T_j$ is a translate of $T_i$. We denote $[T_i]$ the equivalence class
of tile $T_i$, and say $T_i$ is a tile of {\bf type} $[T_i]$.

A {\bf patch} in a tiling is a finite set of its tiles. 
Two patches are said to be {\bf equivalent} if one is a translate of the other,
that is, there is a single translation which takes every tile in one patch
to an equivalent tile in the other patch. The {\bf radius of a patch} is
the radius of the smallest ball containing the patch. 


A tiling is said to have a {\bf finite number of local configurations}, or FLC for short, if there are a finite
number of equivalence classes of patches, up to translation, of any given radius. 

An FLC tiling is {\bf repetitive} if for all $r>0$ there is an $R>0$ such that every patch
of radius $r$ can be found, up to translation, in any ball of radius $R$ in the tiling.
This is equivalent to minimality of the orbit closure of the tiling, see e.g. \cite{Radin},
and was called quasiperiodicity in \cite{Thurston, Ken.thesis}.

In an SAT, the $\phi$-image of each tile type is a well-defined collection
of translates of tile types. If $T_i$ is a tile we can write
$\phi T_i = \cup_j (T_{i_j}+d_{i_j}),$ which is a finite interior-disjoint union.
This subdivision only depends on the type of tile $T_i$,
in the sense that equivalent tiles have equivalent subdivisions.
In particular we let $m_{ij}$ be the number of tiles of type $j$ in the subdivision
of a tile of type $i$. The matrix $\m=(m_{ij})$ is the {\bf subdivision matrix},
it is a nonegative integer matrix which is {\bf primitive}: some power
is strictly positive (by repetitivity of the 
tiling). The leading eigenvalue of $\m$ is the volume expansion
of the SAT, which therefore must be a real Perron number.

Given an SAT, one can select in each of the tile types a point,
called a {\bf control point},
in such a way that the set $\Ck$ of the control points of tiles in a tiling is forward invariant under
$\phi$: $\phi \Ck\subset \Ck$. This can be accomplished as follows \cite{Thurston} (see also \cite[Prop.\,1.3]{Prag}):
for each tile type $[T_i]$, select one tile in its image under expansion and subdivision.
Let the preimage of this tile be $A[T_i]\subset[T_i]$. Then the sequence 
$[T_i],A[T_i],A(A[T_i]),\dots$ nests down to a single point in $[T_i]$, denoted by $c(T_i)$,
which we define to be the control point of $T_i$. For a tile $T = T_i +x$ we let $c(T) = c(T_i)+x$.

\section{Theorem}

The following theorem is stated in \cite{Ken.thesis}.

\begin{theorem} \label{th-main}
Let $\phi$ be a diagonalizable (over $\Compl$) expanding linear map on $\R^n$,
and let $\Tk$ be a self-affine tiling of $\R^n$ with expansion $\phi$. Then

{\bf (i)} every eigenvalue of $\phi$ is an algebraic integer;

{\bf (ii)}
 if $\lam$ is an eigenvalue of $\phi$ of multiplicity $k$ and $\gam$ is
an algebraic conjugate of $\lam$, then either $|\gam| < |\lam|$, or $\gam$
is also an eigenvalue of $\phi$ of multiplicity greater than or equal to $k$.
\end{theorem}

The proof is based on the arguments of Thurston \cite{Thurston} and Kenyon
\cite{Ken.thesis}, but we fill several gaps in those arguments
and provide a great deal more detail.
In particular, Lemmas~\ref{lem-misha} and \ref{lem-Jensen} have no analogs in
\cite{Thurston,Ken.thesis}. It should be pointed out that the corresponding parts of \cite{Thurston} and \cite{Ken.thesis}
have never appeared in refereed publications, but have been widely cited and used in the literature on tilings and tiling dynamical systems.

\medskip

By appropriate choice of a basis, we can assume that  the linear map $\phi$ has the real canonical form, see \cite[Th.\,6.4.2]{HS}.
Since $\phi$ is diagonalizable over $\C$, this means that we have a direct sum decomposition
\begin{equation} \label{bigoplus}
\R^n = \bigoplus_{i=1}^p E_i
\end{equation}
 into invariant subspaces
associated with eigenvalues $\lam_i$ of $\phi$, where we count eigenvalues, having non-negative imaginary part, with multiplicities.  
For a real eigenvalue $\lam_i$, the subspace $E_i$ is one-dimensional, 
and $\phi|_{E_i}$ acts
as multiplication by $\lam_i$. For a non-real eigenvalue $\lam_i$,  the subspace $E_i$ is two-dimensional. Identifying it with a complex plane, we
get that $\phi|_{E_i}$ acts as  multiplication by the complex number $\lam_i$, in other words, as a composition of a dilation and a rotation.
We can define a norm on $\|\cdot\|$ on $\R^n$ such that
\begin{equation} \label{def-norm}
\|x\| = \max_i\|x_i\|\ \ \mbox{for}\ 
x= \sum_{i=1}^p x_i,\ x_i \in E_i,\ \ \|\phi x_i\| = |\lam_i| \|x_i\|
\end{equation}
(here $\|x_i\|$ is just the Euclidean norm on $E_i$ in our basis).

\medskip

{\em Beginning of the proof.} 
Let $\Ck = \Ck(\Tk)$ be a set of control points of the tiling $\Tk$.
Recall that $\phi(\Ck)\subset \Ck$ by construction. Consider $J = \langle \Ck
\rangle$, the free Abelian group generated by $\Ck$. 
It is easy to see that $J$ is finitely generated. Indeed, let
\begin{equation} \label{def-Psi}
\Psi:= \{c(T')-c(T):\ T,T'\in \Tk,\ T\ne T',\ T\cap T'\ne \es\}.
\end{equation}
The set $\Psi$ is finite by FLC, and $J$ is generated by $\Psi$ and
an arbitrary control point
(we can get from it to any control point by moving ``from
neighbor to neighbor'').
Let us fix free generators $v_1,\ldots,v_N$ of $J$. These are 
vectors in $\R^n$; of course, they need not be in $\Ck$. They span $\R^n$,
since $\Ck$ is relatively dense. Note that the choice of the generators is non-unique; in fact, we will need to choose them in a specific way at the end of the proof.
However, for now any generators will do. Let $V$ be the matrix
$V = [v_1\ldots v_N]$. This is a $n\times N$ matrix of rank $n$. By the
definition of free generators,
for every $\xi\in J$ there exists a unique $a(\xi) \in \Z^N$ such that
\begin{equation} \label{eq-address}
\xi = Va(\xi).
\end{equation}
We call $\xi\mapsto a(\xi)$ the ``address map.''
Observe that
\begin{equation} \label{eq-span}
\Span_{\R}\{a(\xi):\ \xi\in \Ck\} = \R^N.
\end{equation}
Indeed, $J$ is generated by $\Ck$, hence every $v_j$ is an integral
linear combination of control points, and $a(v_j)$ is the $j$th unit vector
in $\R^N$.

\begin{lemma} \label{lem-lip} The address map is uniformly
Lipschitz on $\Ck$:
there exists $L_1>0$ such that
\begin{equation} \label{eq-lip1}
\|a(\xi)-a(\xi')\| \le L_1 \|\xi - \xi'\|\ \ \mbox{for all}\ \xi,\xi'\in
\Ck.
\end{equation}
\end{lemma} 

This lemma is a special case of the implication (i) $\Rightarrow$ (v) in  \cite[Th.\,2.2]{Lag1}.
Note that the address map is usually
not even continuous on $J$, since $J$ is not discrete in $\R^n$ 
unless we have a 
``lattice tiling,''
whereas the range of the address map is a subset of the integer lattice in
$\R^N$.

\medskip

Observe that
$\phi \Ck\subset \Ck$ implies $\phi J \subset J$, hence there exists
an integer $N\times N$ matrix $M$ such that
\begin{equation} \label{eq1}
\phi V = V M.
\end{equation}
In other words, we have the commutative diagram (where $i$ indicates the natural inclusion)
$$\begin{CD}
\Z^N @>i>> \R^N @>M>> \R^N @<i<< \Z^N \\
@AA\mbox{$\displaystyle{a}$}A            @V\mbox{$\scriptstyle{V}$}VV       @VV\mbox{$\scriptstyle{V}$}V            @AA\mbox{$\displaystyle{a}$}A \\
J @>i>> \R^n @>\phi>> \R^n @<i<< J
\end{CD}
$$ 

\smallskip

\noindent
For every (complex) eigenvalue $\lam$ of $\phi$
we can find a (complex) left eigenvector $e_\lam$ of $\phi$ corresponding to
$\lam$. Then
$e_\lam V$ is a left eigenvector for $M$ corresponding to $\lam$
(note that $e_\lam V\ne 0$ since $V$ has maximal possible rank $n$).
This proves (i): every eigenvalue of $\phi$ is also an eigenvalue of $M$,
hence an algebraic integer.
Note also that (\ref{eq1}) implies
\begin{equation} \label{eq2}
a(\phi \xi) = Ma(\xi),\ \ \forall\ \xi \in J.
\end{equation}

\begin{lemma} \label{lem-diag}
The matrix $M$ is diagonalizable over $\Compl$.
\end{lemma}

{\em Proof.} 
Recall that $J$ is a free $\Z$-module, on which $\phi$ acts as an endomorphism, and
$M$ is the matrix if this endomorphism in the basis $\Vk:=\{v_1,\ldots,v_N\}$. Note
that $\Q\cdot J$ is a vector space over $\Q$, and $\Vk$ is also a basis of this
vector space. Then $\phi$ induces a linear transformation of $\Q\cdot J$, whose matrix in the
basis $\Vk$ is also $M$. 

Consider the decomposition (\ref{bigoplus}) of $\R^n$ into real eigenspaces
$E_i$ corresponding to the eigenvalues $\lam_i$ of $\phi$.
Decomposing the vectors $v_j$ (the generators of $J$)
in terms of $E_i$ yields
$$
J \subset J':= \bigoplus_{i=1}^p J_i e_i,
$$
where $e_i\in E_i$ and $J_i$ is a finitely-generated $\Z[\lam_i]$-module.
(Here we identify two-dimensional subspaces
$E_i$ with a complex plane on which $\phi$ acts
as multiplication by $\lam_i$.)
Then $\Q\cdot J_i$ is a vector space over $\Q$ and over $\Q(\lam_i)$ (a field).
Let $\{y^{(i)}_1,\ldots,y^{(i)}_{r_i}\}$ be a basis of $\Q\cdot J_i$ over $\Q(\lam_i)$.
Let $n_i$ be the degree of the algebraic integer $\lam_i$. Then 
$\{\lam_i^s y_k^{(i)}:\ 0\le s \le n_i-1,\ 1\le k \le r_i,\ i\le p\}$ is a basis for 
the vector space $\Q\cdot J'$ over $\Q$.
In this basis,  the linear transformation induced by $\phi$ 
has a block matrix, whose every block is
a companion matrix of the minimal polynomial of one of the $\lam_i$'s.
This matrix is diagonalizable over $\Compl$, since the minimal polynomial
has no repeated roots. Finally, we note that the linear transformation induced by
$\phi$ on $\Q\cdot J$ is a restriction of the one which is induced on $\Q\cdot J'$, hence
its matrix, $M$, is diagonalizable as well. \qed

\medskip

Now suppose that $\gam$ is a conjugate of $\lam$, $\gam \ne \lam,\ov{\lam}$,
 and $|\gam|>1$.
Then $\gam$ is an
eigenvalue of $M$. Let $U_\gam$ be the (real) eigenspace for $M$
corresponding to $\gam$.
By Lemma~\ref{lem-diag}, there is a projection $\pi_\gam$ from
$\R^N$ to $U_\gam$ commuting with $M$. By definition, the only
eigenvalues of $M|_{U_\gam}$ are $\gam$ and $\ov{\gam}$ (if $\gam$ is nonreal).
Thus, we can fix a norm on $U_\gam$ satisfying
\begin{equation} \label{ganorm}
\|My\| = |\gam|\,\|y\|,\ \ \ y\in U_\gam.
\end{equation}
Consider the
mapping $f_\gam:\,\Ck\to U_\gam$ given by
\begin{equation} \label{def-fgam}
f_\gam(\xi) = \pi_\gam a(\xi),\ \ \ \xi \in \Ck.
\end{equation}
We would like to extend $f_\gam$ to the entire space $\R^n$. We let
\begin{equation} \label{def-fgam2}
f_\gam(\phi^{-k}\xi) = M^{-k} f_\gam(\xi),\ \ \ \xi \in \Ck.
\end{equation}
This is well-defined since $M$ is invertible on $U_\gam$, and unambiguous
by (\ref{eq2}), since
$\pi_\gam M = M \pi_\gam$. This way we have $f_\gam$ defined on a dense
set
$$
\Ck_\infty:=\bigcup_{k=0}^\infty \phi^{-k} \Ck.
$$

Our goal is to show that $f_\gam$ is uniformly continuous on
$\Ck_\infty$, hence can be extended to all of $\R^n$. In fact, it is
H\"older-continuous. Let $\lam_{\max}$ be the eigenvalue of $\phi$ of 
maximal modulus.
We use the norm (\ref{def-norm}) on $\R^n$.
Denote $B_r(x) = \{y\in \R^n:\ \|y-x\|<r\}$ and let $B_r := B_r(0)$.

\begin{lemma} \label{lem-hoeld} The map $f_\gam$ is 
H\"older-continuous on $\Ck_\infty$: there exists $r>0$ and $L_2>0$ such that 
for any $\xi_1,\xi_2 \in \Ck_\infty$, with $|\xi_1-\xi_2|<r$ we have
\begin{equation} \label{eq-hold}
\|f_\gam(\xi_1) - f_\gam(\xi_2)\| \le L_2 \|\xi_1 - \xi_2\|^\alpha,\ \ \ 
\mbox{for}\ \alpha = \frac{\log|\gam|}{\log|\lam_{\max}|}\,.
\end{equation}
\end{lemma}

{\em Proof.} 
Let $r>0$ be such that for every $x\in \R^n$ the ball $B_r(x)$ is covered by
a tile containing $x$ and its immediate neighbors; this is possible by FLC.
Assume that $\delta = \|\xi_1-\xi_2\|<r$ and 
$\xi_i = \phi^{-k} c_i$ for some $c_i \in \Ck$ and $k\in \Nat$.
Define $\ell$ to be the smallest
positive integer such that
$$
\phi^k B_\delta(\phi^{-k} c_1) \subset \phi^\ell B_r(\phi^{-\ell}c_1).
$$
Since $\ell \le k$, the last inclusion
is equivalent to $|\lam_{\max}|^{k-l}\delta \le r$, so we have
\begin{equation} \label{eq-es1}
|\lam_{\max}|^{-1}(r/\delta) \le |\lam_{\max}|^{k-\ell} \le r/\delta.
\end{equation}
Observe that
$$
c_2 \in \phi^k \ov{B_\delta}(\phi^{-k} c_1) \subset 
\phi^\ell \ov{B_r}(\phi^{-\ell} c_1),
$$
so $\phi^{-\ell}c_1$ and $\phi^{-\ell}c_2$ are in the same or in the
neighboring tiles of $\Tk$ by the choice of $r$.
We claim that there exists a finite set $W\subset J$, independent of 
$c_1,c_2$, such that
\begin{equation} \label{eq-repr}
c_2 - c_1 = \sum_{i=0}^\ell \phi^i w_i
\end{equation}
for some $w_i\in W$ (of course, $w_i$, as well as $\ell$,
depend on $c_1,c_2$).
This is standard, but we provide a proof for completeness.

Let $T_i\in \Tk$ be such that $c_i = c(T_i)$, $i=1,2$. By the 
definition of SAT, there is a (unique) tile $T_i^{(1)}\in \Tk$ such that
$\phi T_i^{(1)}\supset T_i^{(0)} := T_i$. Iterating this, we obtain a sequence
of $\Tk$-tiles $T_i^{(j)}$, for $j\ge 0$, such that $\phi T_i^{(j)} \supset
T_i^{(j-1)}$, for $j\ge 1$ and $i=1,2$.
Note that $T_i^{(\ell)} \supset \phi^{-\ell} T_i^{(0)} \ni \phi^{-\ell} c_i$,
hence $T_1^{(\ell)}$ and $T_2^{(\ell)}$ either coincide or are adjacent.
We have
\begin{eqnarray*}
c_2-c_1  & = & \sum_{j=0}^{\ell-1} \left[\left( \phi^j c(T_2^{(j)}) - 
\phi^{j+1} c(T_2^{(j+1)})\right)-
\left( \phi^j c(T_1^{(j)}) - \phi^{j+1} c(T_2^{(j+1)})\right) \right] \\
& + & \phi^\ell c(T_2^{(\ell)}) - \phi^\ell c(T_1^{(\ell)}). 
\end{eqnarray*}
This implies (\ref{eq-repr}), since the set
$$
\{c(T') - \phi c(T''):\ T',T''\in \Tk,\ T' \subset \phi T''\}
$$
is finite by FLC, as well as the set $\Psi$ from (\ref{def-Psi}), to which
$w_\ell$ belongs.

Now we can write, using (\ref{def-Psi}),
the additivity of the address map on $J$, and (\ref{eq2}),
\begin{eqnarray*}
f_\gam(c_1) - f_\gam(c_2) & = & \pi_\gam a(c_2-c_1) \\
& = & \pi_\gam a\left(\sum_{i=0}^\ell \phi^i w_i \right) \\
& = & \sum_{i=0}^\ell M^i \pi_\gam a(w_i).
\end{eqnarray*}
Thus, in view of (\ref{def-fgam2}) and (\ref{ganorm}),
\begin{eqnarray*}
\|f_\gam(\phi^{-k} c_2) - f_\gam(\phi^{-k} c_1)\| & = & 
\|M^{-k} (f_\gam(c_1) - f_\gam(c_2))\| \\
& = & |\gam|^{-k} \|f_\gam(c_1) - f_\gam(c_2)\| \\
& = & |\gam|^{-k} \left\|\sum_{i=0}^\ell M^i \pi_\gam a(w_i)\right\| \\
& \le & |\gam|^{-k} \sum_{i=0}^\ell |\gam|^i \|\pi_\gam a(w_i)\| \le L'
|\gam|^{\ell-k},
\end{eqnarray*}
where $L' = \frac{|\gam|}{|\gam|-1} \max_{w\in W} \|a(w)\|$.
In view of (\ref{eq-es1}),
$$
|\gam|^{\ell-k}= (|\lam_{\max}|^{\ell-k})^\alpha \le 
(|\lam_{\max}|\delta/r)^\alpha
= \const\cdot \|\xi_1 - \xi_2\|^\alpha,
$$
so we obtain the desired inequality. \qed

\medskip

Now we extend $f_\gam$ by continuity and obtain a function
$f_\gam:\, \R^n\to U_\gam$. 
Observe that
\begin{equation} \label{eq-conj}
f_\gam \circ \phi = M \circ f_\gam,
\end{equation}
since this holds on the dense set $\Ck_\infty$. We also have the following
property. 

\begin{lemma} \label{lem-lip2}
Let $E_\th$ be the real invariant subspace of $\phi$
corresponding to an eigenvalue
$\th$ and suppose that $|\gam| \ge |\th|$. Then $f_\gam|_{{E_\th}+x}$ is
Lipschitz for any $x\in \R^n$, with a uniform constant $2L_1$ (where 
$L_1$ is the constant in Lemma \ref{lem-lip}).
If $|\gam| > |\th|$, then
$f_\gam|_{{E_\th}+x}$ is
constant for any $x\in \R^n$.
\end{lemma}

{\em Proof.} Let $\xi_1,\xi_2\in \R^n$ be such that $\xi_2-\xi_1 \in E_\th$.
By (\ref{eq-conj}), we have for $k\in \Nat$,
\begin{eqnarray*}
\|f_\gam(\xi_1) - f_\gam(\xi_2)\| & = & \|M^{-k} (f_\gam(\phi^k \xi_1)
- f_\gam(\phi^k \xi_2))\| \\
& = & |\gam|^{-k} \|f_\gam(\phi^k \xi_1) - f_\gam(\phi^k \xi_2)\|.
\end{eqnarray*}
Let $c_i$ be a nearest control point to $\phi^k\xi_i$; its distance to
$\phi^k\xi_i$ is at most $d_{\max} = \max\{\diam(T):\ T\in \Tk\}$.
If $k$ is so large that $\|\phi^k \xi_1 - \phi^k \xi_2\|>2d_{\max}$, then
$\|c_1-c_2\| < 2 \|\phi^k \xi_1 - \phi^k \xi_2\|$, and we have by
uniform continuity of $f_\gam$, Lemma~\ref{lem-lip}, and (\ref{def-norm}),
with a uniform
constant $C_3$: 
\begin{eqnarray*}
\|f_\gam(\phi^k \xi_1) - f_\gam(\phi^k \xi_2)\| & \le & 
C_3 + \|f(c_1) - f(c_2)\| \\
& \le & C_3 + L_1\|c_1-c_2\| \\
& \le & C_3 + 2
L_1\|\phi^k \xi_1 - \phi^k \xi_2\| \\
& = & C_3 + 2L_1 |\th|^k \|\xi_1 - \xi_2\|.
\end{eqnarray*}
Thus, 
$$
\|f_\gam(\xi_1) - f_\gam(\xi_2)\| \le C_3|\gam|^{-k} + 2L_1 (|\th|/|\gam|)^k
\|\xi_1 - \xi_2\|.
$$
The lemma follows by letting $k\to\infty$. (Recall that $|\gam| \ge |\th|> 1$.)
\qed

\begin{lemma} \label{lem-loc}
The function $f_\gam$ depends only on the tile type in $\Tk$ up to an
additive constant:
if $T, T+x \in \Tk$ and $\xi \in T$, then
\begin{equation} \label{eq-loc}
f_\gam(\xi+x) = f_\gam(\xi) + \pi_\gam a(x).
\end{equation}
\end{lemma}

Observe that $x\in \Ck - \Ck$, so $a(x)$ is defined, but we cannot write
$\pi_\gam a(x) = f_\gam(x)$, since we do not necessarily have $x\in \Ck$.

\medskip

{\em Proof.}
It is enough to check (\ref{eq-loc}) on a dense set. Suppose
$\xi = \phi^{-k} c(S) \in T$ for some $S\in \Tk$. Then $S\subset \phi^k T$
and
$S + \phi^k x \subset \phi^k(T+x)$ so $S + \phi^k x\in \Tk$. Thus,
\begin{eqnarray*}
f_\gam(\xi+x) & = & f_\gam(\phi^{-k} c(S) + x) \\
& = & f_\gam (\phi^{-k} c(S + \phi^kx)) \\
& = & M^{-k} f_\gam(c(S + \phi^kx)) \\
& = & M^{-k} f_\gam(c(S)) + M^{-k} \pi_\gam a(\phi^k x) \\
& = & f_\gam(\xi) + \pi_\gam a(x),
\end{eqnarray*}
as desired. Here we used the definition of $f_\gam$ on $\Ck$ and
(\ref{eq2}). \qed

\begin{lemma} \label{lem-misha}
If $|\gamma|\geq|\lam|$ then $f_\gam|_{{E_\lam}+x}$ is a constant function
for any $x\in \R^n$.
\end{lemma}

{\em Proof.}
By Lemma~\ref{lem-lip2}, this holds if $|\gam|>|\lam|$, so it remains to
consider the case $|\gam|=|\lam|$. We know that for 
all $x\in \R^n$, the restriction $f_\gam|_{{E_\lam}+x}$ is Lipschitz, hence 
a.e. differentiable by Rademacher's Theorem.
It follows that 
$$
D(x)u := \lim_{t\to 0} \frac{f_\gam(x+tu)-f_\gam(x)}{t}
$$
exists for a.e.\ $x\in \R^n$ for all $u\in E_\lam$, and is a linear
transformation in $u$ (from $E_\lam$ to $U_\gam$).
Moreover, $D(x)$ is measurable in $x$, since it is a limit of continuous
functions.
Since $D(x)$ is the total derivative, we have
\begin{equation}\label{idi2}
\lim_{k\to\infty}\left( \sup_{u\in E_\lam,\ 0< \|u\| < 1/k}
\frac{\|f_\gam(x+u)-f_\gam(x) - D(x)u\|}{\|u\|}\right) =0\ \ \ \mbox{for a.e.}\ x\in \R^n.
\end{equation}
The functions in parentheses are measurable and converge a.e., hence by
Egorov's Theorem they converge uniformly on a set of positive measure.
Uniform convergence means that there exists  a sequence of positive integers
$N_k \uparrow \infty$ such that 
\begin{eqnarray*}
\Omega & := & \{\xi\in \R^n:\ 
\|f_\gam(\xi+u)-f_\gam(\xi) - D(\xi)u\| \le \|u\|/k\\
& & \forall\,u \in B_{1/N_k} \cap E_\lam,\ \mbox{for all $k$} \}
\end{eqnarray*}
has positive Lebesgue measure.
We claim that $\Omega$ has full Lebesgue measure.

Observe that if $T,T+x\in \Tk$ and $\xi\in T^\circ$, then
\begin{equation} \label{idi3}
\xi\in \Omega\ \Rightarrow\ \xi+x\in \Omega
\end{equation}
by Lemma~\ref{lem-loc}.
Furthermore, by (\ref{eq-conj}) we have
$
D(\phi\xi) = MD(\xi) \phi^{-1}
$
and, denoting $v=\phi u$, for all $v\in B_{|\lam|/N_k}\cap E_\lam$,
\begin{eqnarray*}
\|f_\gam(\phi \xi + v) - f_\gam(\phi \xi) - D(\phi \xi) v\| & = & 
\|M(f_\gam(\xi + u) - f_\gam(\xi)) - D(\xi)u) \|\\
& = & |\gam|\cdot \|f_\gam(\xi+u)-f_\gam(\xi) - D(\xi)u\| \\
& \le & |\gam|\cdot \|u\|/k = |\lam|\cdot \|u\|/k = \|v\|/k,
\end{eqnarray*}
where we used that $\phi|_{E_\lam}$ expands the norm by a factor of $|\lam|$.
This shows that $\phi(\Omega) \subset \Omega$.

We will need a version of Lebesgue-Vitali 
Density Theorem where the differentiation
basis is not the set of balls but rather the collection of sets of the
form $\phi^{-k}B_1$, $k\ge 0$, and their translates.
It is a well-known fact in Harmonic Analysis that such sets form a density
basis, for any expanding linear map $\phi$ (even non-diagonalizable), 
see \cite[pp.\,8-13]{Stein} or \cite[pp.\,11-14]{PughShub}.
Let $y$ be a density point of $\Omega$, 
i.e., denoting the Lebesgue measure by $m$,
$$
m(\Omega \cap \phi^{-k} B_1(\phi^k y)) \ge (1-\eps_k) m(\phi^{-k} B_1)\ \ \ 
\mbox{for some}\ \eps_k\to 0.
$$
Denote by $[B_1(x)]^\Tk$ the patch consisting of those tiles which intersect
$B_1(x)$. By repetitivity, there exists $R>0$ such that
$B_R$ contains a translate of  $[B_1(x)]^\Tk$ for
every $x\in \R^n$. 
Let $y_k \in B_R$ be such that $[B_1(y_k)]^\Tk$ is a translate of
$[B_1(\phi^k y)]^\Tk$.
Then
\begin{eqnarray*}
m(\Omega\cap B_1(y_k)) & = & m(\Omega \cap B_1(\phi^k y)) \\
               & \ge & m(\phi^k \Omega \cap B_1(\phi^k y)) \\
               & = & |\det \phi|^k  m(\Omega \cap \phi^{-k} B_1(\phi^k y)) \\
                & \ge & |\det \phi|^k (1-\eps_k) m(\phi^{-k} B_1)
               =(1-\eps_k) m(B_1).
\end{eqnarray*}
We used (\ref{idi3}) and  $\phi^k \Omega \subset \Omega$ in the first two
displayed lines above.
Let $y'$ be a limit point of $y_k$. Then we have $m(\Omega\cap B_1(y')) = 
m(B_1)$. Thus, $\Omega$ is a set of full measure in $B_1(y')$, and by
expansion and translation we conclude that $\Omega$ has full measure in $\R^n$,
completing the proof of the claim.

Now choose $\ell_k$ so that $|\lam|^{\ell_k} > N_k$. We have
\begin{eqnarray*}
\zeta \in \phi^{\ell_k} \Omega & \Rightarrow &
                 \|f_\gam(\zeta+v) - f_\gam(\zeta) - D(\zeta)
                 v\| \le \|v\|/k\\
                 & & \mbox{for all}\ v \in \phi^{\ell_k} (B_{1/N_k} \cap
              E_\lam) \supset B_1 \cap E_\lam.
\end{eqnarray*}
We know that $\Omega'=\bigcap_{k\ge 1} \phi^{\ell_k} \Omega$ 
has full measure, hence it
is dense. For any $\xi \in \R^n$ choose a sequence $\xi_k\to \xi$
such that $D(\xi_k)$ converges (this is possible since $\|D(\xi)\|\le 2L_1$
by Lemma~\ref{lem-lip2}).
Passing to
the limit, we obtain that
$$
f_\gam(\xi + v) = f_\gam(\xi) + D(\xi) v,\ \ \ 
\mbox{for all}\ v \in B_1\cap E_\lam.
$$
This shows 
that $f$ is affine linear on every $E_\lam$ slice:
$$
f_\gam(\xi + v) = f_\gam(\xi) + D(\xi) v,\ \ \ \mbox{for all}\ v \in E_\lam,
$$
and $D(\xi) = D(\xi')$ whenever $\xi'-\xi \in E_\lam$.
Taking $\xi=0$ we see that $f_\gam|_{E_\lam}$ is linear. It intertwines
$\phi|_{E_\lam}$ and $M|_{U_\gam}$. But 
$\{\gam,\ov{\gam}\}\cap \{\lam,\ov{\lam}\}=\es$ which are the
eigenvalues of $\phi|_{E_\lam}$ and $M|_{U_\gam}$ respectively, hence
the only possibility is $f_\gam|_{E_\lam} \equiv 0$. Since $f_\gam$ 
is uniformly continuous
on $\R^n$ and $f_\gam|_{x+E_\lam}$ is affine linear, we obtain that 
$f_\gam|_{x+E_\lam}\equiv \const(x)$. \qed

\medskip

To motivate the conclusion of the proof, we start with a heuristic discussion.
Assume that $|\gam|\ge |\lam|$ for the rest of the proof.
So far, we have proved that $f_\gam$ is affine linear on the slices $x+E_\lam$. Suppose we could show that $f_\gam$ is linear on $\R^n$. Then we could conclude as follows:
$f_\gam \circ \phi = M \circ f_\gam$ and $f_\gam(\R^n) = U_\gam$ (the latter follows from (\ref{eq-span}) and the definition of $f_\gam$) would imply that
$\phi$ restricted to a linear subspace and $M|_{U_\gam}$ are linearly conjugate:
$$
\begin{CD} 
U_\gam \subset \R^N @>M>> U_\gam \subset \R^N \\
@A\mbox{$f_\gam$}AA                   @A\mbox{$f_\gam$}AA \\
\R^n @>\phi>> \R^n
\end{CD}
$$
and hence $\gam$ is an eigenvalue of $\phi$ of multiplicity at least $\dim\,U_\gam \ge \dim\,E_\lam$, as desired.

This scheme does work, but with some modifications. We are able to show that
$f_\gam$ is affine linear in some, but possibly not all, directions
complementary to $E_\lam$. It is linear in directions for which the differences
between control points for tiles of the same type project densely. 

Let $\Xi = \Xi(\Tk)$ denote the set of translation vectors between
tiles of the same type and let $P_\lam$ be the projection from $\R^n$ to 
$E_\lam$ commuting with $\phi$ (note that the projection $\pi_\gam$ acts
in another space, $\R^N$).

Consider the set $(I-P_\lam)\Xi$,
that is, the projection of $\Xi$ onto the other eigenspaces of $\phi$.
This projection may look like a lattice
in some directions and fail to be discrete in other
directions. We consider the directions in which this set is not discrete;
more precisely, those directions
in which there are arbitrarily small nonzero vectors in
$(I-P_\lam)\Xi$, and denote the span of these directions $E'$. What we will
prove is that $f_\gam$ is affine linear on all $E'$ slices,
and hence all $E'\oplus E_\lambda$ slices. We will then show that the subspace 
$E:=E'\oplus E_\lam$ is $\phi$-invariant and is spanned by the vectors of $\Xi$ contained
in it. This will allow us to essentially restrict the entire construction to $\phi|_E$ and conclude as indicated above, using that $f_\gam|_E$ is linear.

Now let us be more formal
and for each $\eps>0$ define $E_\eps\subset\R^n$ to be the subspace
$$E_\eps=\Span_{\R}(B_\eps\cap(I-P_\lam) \Xi)\subset E_\lambda^\perp\subset\R^n,$$
where $E_\lambda^\perp$ is the $\phi$-invariant subspace complementary to $E_\lam$.
Further, consider
$$
E': = \bigcap_{\eps>0} E_\eps.
$$
We have $\phi \Xi \subset \Xi$ and $P_\lam \phi = \phi P_\lam$, hence
$$
\phi((I-P_\lam) \Xi) \subset (I-P_\lam) \Xi.
$$
Note that $E_\eps$ are decreasing linear subspaces of 
$E_\lambda^\perp\subset\R^n$, hence 
$E'= E_\eps$ for some $\eps>0$, and so $E' = E_{\eps'}$ for all
$0 < \eps' \le \eps$. Since $\phi E_{\eps'} \subset E_{c\eps'}$ for
$c=\|\phi\|$ we see that $E'$ is $\phi$-invariant.
We then define
$$
E:= E'+E_\lam.
$$

\begin{lemma} \label{lem-Jensen}
$f_\gam|_{E+x}$ is affine linear for every $x\in \R^n$.
\end{lemma}

{\em Proof.}
Choose $\eps$ so that $E' = E_\eps$.
Let $\eps'<\eps$ and define
$$E'':= \Span(B_{\eps'}\cap (I-P_\lam)(\Ck_1-\Ck_1))$$
where $\Ck_1$ is the set of control points of tiles of type 1 (of course, we could equally well choose another tile type).
First we claim that
\begin{equation} \label{eq-claim1}
E' = E''.
\end{equation}
Indeed, $\Ck_1-\Ck_1 \subset \Xi$ hence $E''\subset E'$.
Choose $\ell$ so large that $\phi^\ell \Xi \subset \Ck_1-\Ck_1$; such an
$\ell$ exists by primitivity of the tile substitution
(the $\ell$-th power of the substitution
of any tile contains tiles of all types).
We then have
$$
E' = \phi^\ell E' = \phi^\ell E_{\eps'/\|\phi\|^\ell} \subset
 \Span(B_{\eps'}\cap (I-P_\lam)\phi^\ell \Xi) \subset E''.
$$
The claim is proved.

Now suppose $x\in \Ck_1-\Ck_1$, so
there exists $T\in \Tk$ of type 1 such that $T+x\in \Tk$.
By Lemma~\ref{lem-loc},
$$
\xi\in T \ \Rightarrow\ f_\gam(\xi+x) = f_\gam(\xi) + \pi_\gam a(x).
$$
But Lemma~\ref{lem-misha} implies that $f_\gam(\xi+x) = f_\gam(\xi + x-P_\lam x)$, so
\begin{equation} \label{idi4}
f_\gam(\xi + (I-P_\lam)x) = f_\gam(\xi) + \pi_\gam a(x)\ \ \ \mbox{for}\ \xi \in T.
\end{equation}

We want to show that $f_\gam$ is affine linear on all $E$-slices. Since $f_\gam$ is
constant on all $E_\lam$-slices by Lemma~\ref{lem-misha}, it is enough to verify that $f_\gam$ is affine
linear on all $E'$-slices (recall that $E = E'+E_\lam$).
Fix a small $\eps'$ as in (\ref{eq-claim1}) and select a basis of $E'$ of the
form $y_i = (I-P_\lam) x_i\in B_{\eps'}$, with $x_i \in \Ck_1-\Ck_1$,
for $i=1,\ldots, \dim E'$. Now for any $\xi$ in the interior of $T$, such that
$B_r(\xi) \subset T$, we
obtain from (\ref{idi4}):
$$
f_\gam\Bigl(\xi + \sum_i b_i y_i\Bigr) = f_\gam(\xi) + \sum_i b_i \pi_\gam a(x_i),
$$
for all $b_i \in \Z$ such that  $\sum_i b_i y_i\in B_r$.
(Here we should note that, in view of Lemma~\ref{lem-loc}, equality
(\ref{idi4}) transfers to all tiles equivalent to $T$. Since all the $x_i$
are translates between two copies of $T$, we can apply the equality for any
$x_i$ in any of the translates.)
This shows that $f_\gam$ is affine linear on a large
chunk of the lattice in $E'$ generated by small vectors $y_i$,
translated in such a way that $\xi$ becomes the origin.
It is an easy exercise to
pass to the limit as $\eps'\to 0$ and conclude that $f_\gam$ is affine linear
in the $E'$-direction on $B_r(\xi) \cap (E'+\xi)$. To be a bit more precise,
we can verify that
\begin{equation} \label{idi5}
f_\gam\Bigl( \frac{\zeta_1+\zeta_2}{2} \Bigr) = \frac{f_\gam(\zeta_1) + f_\gam(\zeta_2)}{2}
\ \ \ \mbox{for all}\ \zeta_1,\zeta_2 \in B_r(\xi) \cap (E'+\xi).
\end{equation}
Since $f_\gam$ is continuous, this implies that
\begin{equation} \label{idi6}
f_\gam(\zeta) = A_\xi \zeta + b_\xi\ \ \ \mbox{for all}\ \zeta\in B_r(\xi) \cap
(E'+\xi),
\end{equation}
see e.g., \cite[2.1.4]{Aczel},
where it is called the ``Jensen functional equation''. The details are
straightforward.

Since (\ref{idi6}) holds on all slices of $T$, by ``expanding and
translating'' with the help of (\ref{eq-conj}) and Lemma~\ref{lem-loc},
we obtain the claim of the lemma. \qed

\medskip

\begin{lemma} \label{lem-KS3}
$$
E=\Span_\R((\Ck-\Ck) \cap E).
$$
\end{lemma}

{\em Proof.} Denote $W:= \Span_\R((\Ck-\Ck) \cap E).$ First we show that $E_\lam \subset W$.
Let $w\in E_\lam$. The set $\Ck_1$ (control points of type-1 tiles) is relatively dense in $\R^n$; let $R>0$ be such that
every open ball of radius $R$ hits $\Ck_1$. Let $\xi_j\in \Ck_1$ be such that $\|\xi_j - jw\| < R$ for all
$j\ge 0$. Then
\[ \|(I-P_{\lam})\xi_j\| = \|(I-P_{\lam})(\xi_j-jw)\| \le (1+ \|P_{\lam}\|)R,\ \ j\ge 0.\]
It follows that there exists a sequence of pairs $(i_k,j_k)$, with $i_k-j_k \to +\infty$, such that
\[ \|(I-P_{\lam})(\xi_{i_k} - \xi_{j_k})\| \to 0,\ \ \mbox{as}\ k\to \infty.
\]
Therefore, $(I-P_{\lam})(\xi_{i_k} - \xi_{j_k})\in E'$ for $k$ sufficiently large, and hence $\xi_{i_k} - \xi_{j_k}
\in E$ for $k\ge k_0$. Now,
\[ \| (\xi_{i_k} - \xi_{j_k}) - w(i_k - j_k) \|\le 2R,\]
hence $\zeta_k := (\xi_{i_k} - \xi_{j_k})/(i_k-j_k) \to w$. But $\zeta_k \in W$, for $k\ge k_0$, hence
$w\in W$, since $W$ is closed, being a linear subspace of $\R^n$. 

Now recall that $E'$ is spanned by certain vectors of the form $\xi - P_\lam \xi$, with $\xi \in \Xi \subset \Ck-\Ck$. Since $P_\lam \xi \in E_\lam \subset E$, we have that these
vectors $\xi$ are in $E$, and hence $E' \subset W$. This proves that $E=E'+E_\lam \subset W$, as desired. 
\qed

\medskip

{\em Conclusion of the proof of Theorem~\ref{th-main}.} 
As mentioned earlier, we would like to run the entire construction essentially restricting ourselves to the subspace $E$, which is $\phi$-invariant, contains $E_\lam$, and is spanned
by the vectors of $\Ck-\Ck$ in it. We do not literally do this, because it is not clear what the
intersection of the tiling with $E$ looks like; rather, we make sure that the construction on $\R^n$ is
compatible with this subspace structure.
Recall that at the beginning of the proof we considered the free Abelian group $J = \langle \Ck \rangle$ and its free generators $v_1,\ldots,v_N$.
We will now use a more specific choice of the generators. Namely, let
$$\wtil{J}:= \langle (\Ck-\Ck)\cap E\rangle = \Span_\Z((\Ck-\Ck) \cap E).$$
 Clearly, $\wtil{J}$ is an Abelian subgroup of $J$, 
and $\Span_\R \wtil{J} = E$ by Lemma~\ref{lem-KS3}. Is it possible to choose the free generators for $J$ as an extension of a set of free generators for $\wtil{J}$? Maybe not, but 
we can choose $v_1,\ldots,v_N$, the free generators of $J$, 
so that
$d_1 v_1,\ldots, d_s v_s$ are free generators of $\wtil{J}$
for some positive integers $d_j$ and $s\le N$
(see e.g. \cite[Theorem II.1.6]{hung}). 

Recall that $\phi$ acts on $J$, and on the generators $v_j$ this action is given by an integer matrix $M$. Since $\phi$ also acts on $\wtil{J}$, we claim that
$M = \left( \begin{array}{c|c} \wtil{M} & * \\ \hline 0 & * \end{array} \right)$, where $\wtil{M}$ is an $s\times s$ matrix.
Indeed, $\phi(v_i)$, $i\le N$, is a unique integral linear combination of $\{v_j\}_{j\le N}$, with the coefficients coming from the $i$-th column of $M$.
On the other hand, $\phi(d_i v_i)$, $i\le s$, is an 
integral linear combination of $\{d_jv_j\}_{j\le s}$, since the latter are free generators of $\wtil{J}$.
This implies that $\phi(v_i)$, $i\le s$, is an integral linear combination of $\{d_jv_j\}_{j\le s}$, that is,
\begin{equation} \label{eq-conj2}
\phi [v_1\ldots v_s] = [v_1\ldots v_s]\wtil{M},
\end{equation}
where $\wtil{M}$ is an integral $s\times s$ matrix. Thus, the matrix $M$ is
block upper-triangular, with the upper left corner $\wtil{M}$, as claimed above.

Note that
\begin{equation} \label{eq-span3}
\Span_\R (\{v_j\}_{j\le s}) = \Span_\R (\{d_j v_j\}_{j\le s}) = \Span_\R((\Ck-\Ck) \cap E) = E
\end{equation}
by construction.
By (\ref{eq-conj2}) and (\ref{eq-span3}), there is an $\wtil{M}$-invariant subspace of $\R^s$, on which $\wtil{M}$ acts isomorphically
(linearly conjugate) to $\phi|_E$. Since $E\supset E_\lam$, we obtain that 
$\lam$ is an eigenvalue of $\wtil{M}$, with the multiplicity
greater or equal to $\dim E_\lam$. 
Because $\gam$ is an algebraic conjugate of $\lam$ and $\wtil{M}$ is an integer matrix, we have that $\gam$ is also an eigenvalue
of $\wtil{M}$, with the multiplicity $\ge \dim E_\lam$. 
Let $\wtil{U}_\gam$ be the real invariant subspace of
$\wtil{M}$ corresponding to $\gam$.

Abusing notation a bit, we will identify $\R^s$ with the subspace of $\R^N$ generated by the first $s$ coordinates. Then $\wtil{U}_\gam\subset U_\gam$.

Let $a:\, J \to \Z^N$ be the address map, as in (\ref{eq-address}). Then $a(\wtil{J}) \subset \Z^s$ (using a similar abuse of notation, so that $\Z^s\subset \Z^N$).
By construction,
$$
\Span_\Z\{a(\xi-\xi'):\ \xi,\xi'\in \Ck,\ \xi-\xi'\in E\}=\bigoplus_{j=1}^s d_j \Z \subset \Z^s,
$$
hence
$$
\Span_\R\{a(\xi-\xi'):\ \xi,\xi'\in \Ck,\ \xi-\xi'\in E\}=\R^s.
$$
It follows that 
\begin{equation} \label{eq-duda}
\Span_\R\{\pi_\gam (a(\xi)-a(\xi')):\ \xi,\xi'\in \Ck,\ \xi-\xi'\in E\}=\pi_\gam(\R^s) = \wtil{U}_\gam.
\end{equation}

Recall that $f_\gam:\,\R^n \to \R^N$, defined originally by $f_\gam(\xi) = \pi_\gam (a(\xi))$ on control points,
is uniformly continuous, $f_\gam \circ \phi = M \circ f_\gam$, and 
$f_\gam|_{E+x}$ is affine linear for all $x$ by Lemma~\ref{lem-Jensen}. 
Note that $f_\gam|_E$ is linear, since $f_\gam(0)=0$. 

We claim that $f_\gam(E) \supset \wtil{U}_\gam$. Indeed, every $f_\gam(E+x)$ is a translate of a linear subspace, which must be a  translate of $f_\gam(E)$, by the
uniform continuity of $f_\gam$. It follows that for $\xi,\xi'\in \Ck,\ \xi-\xi'\in E$,
$$
\pi_\gam(a(\xi)-a(\xi')) = f_\gam(\xi)-f_\gam(\xi') \in f_\gam(E),
$$
whence $\wtil{U}_\gam \subset f_\gam(E)$ by (\ref{eq-duda}). The claim is verified.

Since $f_\gam(E)$ contains $\wtil{U}_\gam$, there exists a $\phi$-invariant
subspace $\widetilde{E}\subset E\subset \R^n$, such that $f_\gam$ maps
$\widetilde{E}$ isomorphically onto $\wtil{U}_\gam$:
$$
\begin{CD}
\wtil{U}_\gam \subset \R^s @>\wtil{M}>> \wtil{U}_\gam \subset \R^s \\
@A\mbox{$f_\gam$}AA                   @A\mbox{$f_\gam$}AA \\
\wtil{E} \subset E  @>\phi>> \wtil{E} \subset E
\end{CD}
$$
Thus, the linear map $f_\gam|_{\wtil{E}}$ conjugates $\phi|_{\wtil{E}}$ to 
$\wtil{M}|_{\wtil{U}_\gam}=M|_{\wtil{U}_\gam}$, hence
$\gam$ is an eigenvalue of $\phi$ of multiplicity $\ge \dim E_\lam$,
as desired. \qed

\medskip

{\bf Acknowledgment.} We are grateful to Misha Lyubich for a
suggestion  which helped prove Lemma~\ref{lem-misha}.

\end{document}